# Gradual Variation Analysis for Groundwater Flow

Li Chen
Department of Computer Science and Information Technology
University of the District of Columbia
Email: *lchen@udc.edu*
(202) 274-6301

**Abstract**

Groundwater flow in Washington DC greatly influences the surface water quality in urban areas. The current methods of flow estimation, based on Darcy's Law and the groundwater flow equation, can be described by the diffusion equation (the transient flow) and the Laplace equation (the steady-state flow). The Laplace equation is a simplification of the diffusion equation under the condition that the aquifer has a recharging boundary. The practical way of calculation is to use numerical methods to solve these equations. The most popular system is called MODFLOW, which was developed by USGS. MODFLOW is based on the finite-difference method in rectangular Cartesian coordinates. MODFLOW can be viewed as a "quasi 3D" simulation since it only deals with the vertical average (no *z*-direction derivative). Flow calculations between the 2D horizontal layers use the concept of leakage.

In this project, we have established a mathematical model based on gradually varied functions for groundwater data volume reconstruction. These functions do not rely on the rectangular Cartesian coordinate system. A gradually varied function can be defined in a general graph or network. Gradually varied functions are suitable for arbitrarily shaped aquifers. Two types of models are designed and implemented for real data processing: (1) the gradually varied model for individual (time) groundwater flow data, (2) the gradually varied model for sequential (time) groundwater flow data. In application, we also established a MySQL database to support the related research.

The advantage of the gradually varied fitting and its related method does not need the strictly defined boundary condition as it is required in MODFLOW.

## 1. Introduction

Groundwater flow in DC has greatly influenced the surface water quality in urban areas. The current method of flow estimation mainly uses the ground flow equation, which is a partial differential equation. Software systems such as MODFLOW can only solve 2D equations and pass the data vertically to form a 3D volume. The research on the 3D models of groundwater flow has fundamental and practical importance to hydrogeology. A method used to establish a true 3D Groundwater flow will be very useful to groundwater related research in DC. The new method can also be extended for use in other regions.

A feasible and true 3D model for groundwater flow is essential to groundwater research. This project attempts to establish a 3D model using discrete mathematics, especially graphical and graph-theoretical methods, to compute groundwater flow. It is called the gradual variation method. This method can be combined with the existing technology such as MODFLOW for more accurate results. The research is focused on the groundwater flow of the DC area. We have established a connection to use gradually varied functions in groundwater research. We have accomplished three tasks: (1) extracting of real data from databases of DC areas, (2) storing the data into local database, (3) reconstructing the water-head surfaces for time sequences using gradually varied surface fitting. We have also completed the design of the combined gradually varied fitting using the finite difference method.

The advantage of the gradually varied fitting and its related method does not need the strictly defined boundary condition as it is required in MODFLOW.

The results of the research will provide a reliable source to better understand the groundwater of DC. The images of data flow will indicate the activity of the groundwater flow. Beyond the theoretical achievements, the undergraduate and graduate students of computer science will have the opportunity to learn how to conduct state of the art groundwater research and use various software systems. These students could work in the environmental sciences in the future.

## 2. Problems and Methods

In order to use discrete methods for Groundwater Modeling, we need to solve the following problems: (1) Groundwater flow equations and its discrete forms, (2) Gradually varied functions for groundwater data, (3) Real data preparation, (4) Algorithms, especially fast design, (5) Real data processing and applications.

Much research has already been done to find a discrete model for the groundwater flow equations. The research is mainly based on numerical methods and analytical methods. The finite difference method (FDM) and the finite element method (FEM) are popular in this area. Pruist et al. [22] has indicated that FEM has advantages on local refinement of grid (adaptive mesh generation) due to non-rectangular grids, good accuracy, stability, representation of the spatial variation of anisotropy. However, its computational cost is much larger and relative less manageable in application. In fact, groundwater industry is not like the automobile industry, there is no much need for a good-looking smoothed groundwater level surface.

On the other hand, FDM has simplicity of theory and algorithm, plus easiness of application; however it has the problems of inefficient refinement of grid and poor geometry representation. This is because of the strict use of rectangular grids. In addition, there is no standard method to implement the Neumann boundary condition.

The gradually varied function is supposed to pick up the advantages and overcome the disadvantages in FEM and FDM. So our first task is to investigate the suitableness of gradually varied functions for groundwater data. Then, we must find a connection between the flow equations and gradually varied functions. We also need to design an input data format to store the data in a database.

A gradually varied function is for the discrete system where a high level of smoothness is not a dominant factor. It can be used in any type of decomposition of the domain. It is more flexible than rectangle-cells used in MODFLOW and triangle-cells used in FEFLOW. Because gradual variation does not have strict system requirements, the other mathematical methods and the artificial intelligence methods can be easily incorporated into this method to seek a better solution. Based on the boundary conditions or constraints of the groundwater aquifer, the constraints could be in explicit forms or in differential forms such as the diffusion equations. The gradually varied function exists based on the following theorem: The necessary and sufficient condition of the existance of a gradually varied function is that the change of values in any pair of sample points is smaller or equal to the distance of the points in the pair.

## 3. Background and Related Research

The research of groundwater flow is one of the major topics in groundwater research [1] [2] (this sentence does not make sense and just repeats itself). The groundwater flow equation based on Darcy's Law usually describes the movement of groundwater in a porous medium such as aquifers. It is known in mathematics as the diffusion equation. The Laplace equation (for a steady-state flow) is a simplification of the diffusion equation under the condition that the aquifer has a recharging boundary. The conservation of mass states that for a given increment of time ($\Delta t$) "the difference between the mass flowing in across the boundaries, the mass flowing out across the boundaries, and the sources within the volume, is the change in storage."[2]

$$\frac{\Delta M_{stor}}{\Delta t} = \frac{M_{in}}{\Delta t} - \frac{M_{out}}{\Delta t} - \frac{M_{gen}}{\Delta t} \qquad (1)$$

Its differential form is

$$\frac{\partial h}{\partial t} = \alpha \left[ \frac{\partial^2 h}{\partial x^2} + \frac{\partial^2 h}{\partial y^2} + \frac{\partial^2 h}{\partial z^2} \right] - G$$

(2)

To solve this equation, a grid method is usually used, such as the finite difference or finite element method [3] [4]. Other methods including the analytic element method attempt to solve the equation exactly, but need approximations of the boundary conditions [5][6]; they are mainly used in academic and research labs.

The practical way of calculating is to use numerical methods to solve these equations. The most popular system is called MODFLOW, which was developed by USGS [7]. MODFLOW is based on the finite-difference method on rectangular Cartesian coordinates. MODFLOW can be viewed as a "quasi 3D" simulation since it only deals with the vertical average (no *z*-direction derivative). Flow calculations between the 2D horizontal layers use the concept of leakage.

Gradual variation is a discrete method that can build on any decomposition or networking.  It was originally introduced to solve image processing problems and discrete surface reconstruction [8][9][21]. See Fig 1. When a boundary is known, it can used to fit (solve) for the interal points in any type of linked connection.   The original research

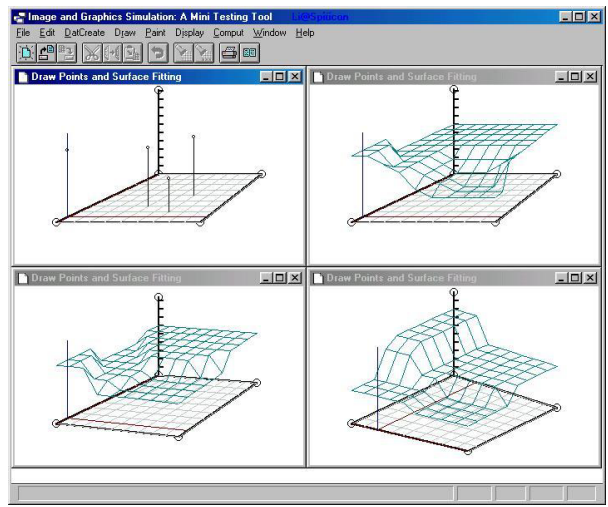

Fig. 1  Examples of gradually varied functions

Gradually Varied Surfaces of this project is to use Darcy's Law or differential constrains to determine the value of the unknown points instead of the random selection of the construction of gradually varied surfaces when there are more than one possible selections [10].  When the determination of the values are uncertain, we will try to use the artificial intelligence methods such as neural networks and genetic algorithms to help us to find near optimal solutions. These types of studies have already been done by many researchers in groundwater flow [11][12].

On the other hand, current research still showed considerable interest in establishing new modeling methods for groundwater flow [13-16]. Both the US and DC governments are concerned with  the  future of groundwater flow research [17][18].  This project on the 3D models of groundwater flow has fundamental and practical importance to hydrogeology. A method used to establish a true 3D model for groundwater flow will be very useful to groundwater related research in DC and in other urban areas.

## 4. Data Preparation

Data Preparation is a very important part of this project, please see appendix or [20] for more information. UDC student Branham has to use PHP to build a web application to access groundwater log data in Virginia and Maryland. Data is stored in MySQL databases.

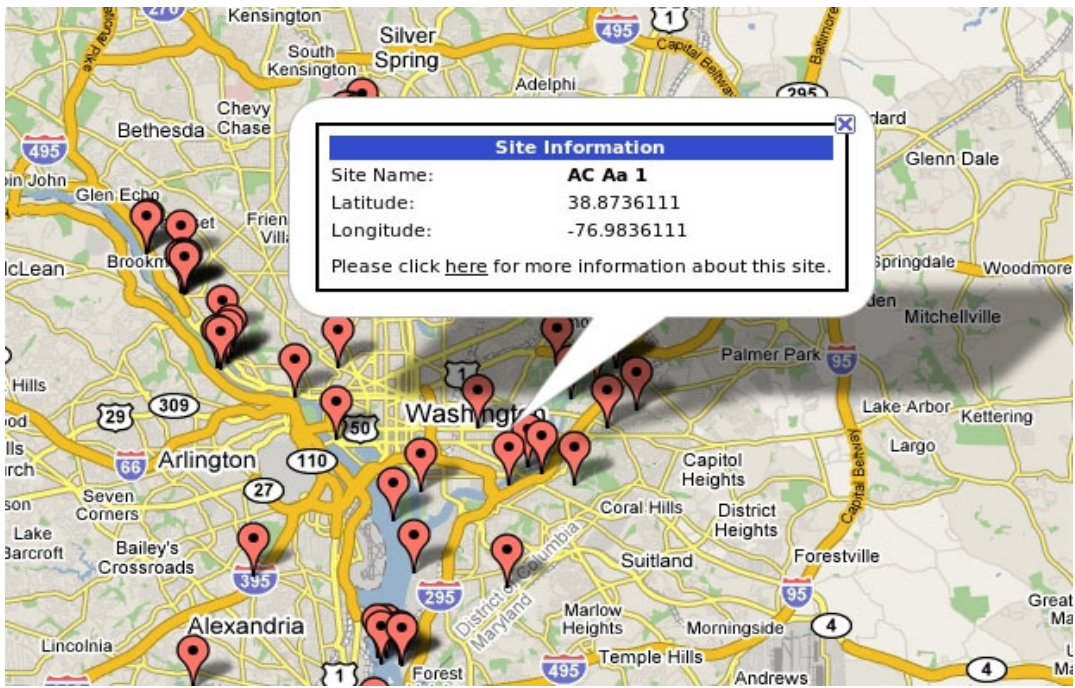

Fig. 2 Web based data collection

The application, in its current form, can be found at:
http://www.travisbranham.com/waterquality/

## 5. Ground Water Surface Fitting

We have done extensive research on data collection and initial data reconstruction using gradually varied functions that are for the discrete system where a high level of smoothness is not a dominant factor. Since they can be used on any type of decomposition of the domain, gradually varied functions are more flexible than rectangle-cells used in MODFLOW and triangle-cells used in FEFLOW. Because gradual variation does not have strict system requirements, the other mathematical methods and the artificial intelligence methods can be easily incorporated into this method to seek a better solution. Based on the boundary conditions or constraints of the groundwater aquifer, the constraints could be in explicit forms or in differential forms such as the diffusion equations. Again we select the gradually varied function because of the following theorem: The necessary and sufficient condition of the existence of a gradually varied function is that the change of values in any pair of sample points is smaller or equal to the distance of the points in the pair.

### 5.1 Individual surface fitting

Gradual variation is a discrete method that can be built on any graph. The gradually varied surface is a special discrete surface. We now introduce this concept.

The Concept of Gradual Variation: Let function f: D$\rightarrow$\{1, 2,…,n\}, if a and b are adjacent in D implies $|f(a)-f(b)| \leq 1$, point (a,f(a)) and (b,f(b)) are said to be gradually varied. A 2D function (surface) is said to be gradually varied if every adjacent pair are gradually varied.

Discrete Surface Fitting: Given J$\subseteq$D, and f: J$\rightarrow$\{1,2,…n\} decide if there is a F: D$\rightarrow$\{1,2,…,n\} such that F is gradually varied where f(x)=F(x), x in J.

**Theorem** (Chen, 1989) The necessary and sufficient conditions for the existence of a gradually varied extension F is: for all x,y in J, $d(x,y) \geq |f(x)-f(y)|$, where d is the distance between x and y in D.

The above theorem can be used for a single surface fitting if the condition in the theorem is satisfied. The problem is that the sample data does not satisfy the condition of fitting. So the original algorithm cannot be used directly for individual surface fitting.

An algorithm (Algorithm A) based on the sample point contribution to the fitting point is created. The core part of the program is given in the following code. This algorithm can produce the surface but do not give good fit in some cases. We have showed the result in the previous report.

```
Algorithm A:
    for (k=0;k<nGuildPoints;k++){
        ii=(latIndex[k]-latMin)/latDet;
        jj=(longtIndex[k]-longtMin)/longtDet;

        distance=sqrt((ii -i)*(ii-i)+(jj -j)*(jj-j));
        temp_j=abs((array[i][j] - dat[k][time]))/Ratio-distance;
        if(temp_j>0){ // not satisfy gvs condition

        if( array[i][j] > dat[k][time])
            temp=-temp_j *Ratio ;
        else
            temp= temp_j *Ratio;
        array[i][j]=array[i][j]+temp ;
    }
```

A new algorithm was developed recently to overcome these problems [23]. In [23], a systematic digital-discrete method for obtaining continuous functions with smoothness to a certain order (C^n) from sample data is designed. This method is also based on gradually varied functions. The new algorithm tries to search for a best solution of the fitting. We have added the component of the classical finite difference method. The result of the new method will be shown in following sections. This method is independent from existing popular methods such as the cubic spline method and the finite element method. The new digital-discrete method has considerable advantages for a large amount of real data applications. This digital method also differs from other classical discrete method that usually uses triangulations. This method can potentially be used to obtain smooth functions such as polynomials through its derivatives $f^{(k)}$ and the solution for partial differential equations such as harmonic and other important equations.

**5.2 Sequential surface fitting and involvement of the flow equation**

The individual surface fitting in the above code is not an adequate method since the relationship of groundwater surfaces at each time is not accounted for in the calculation. In particular, we have to use the flow equation in the sequential surface calculation.

According to flow equation (2),
$$h2-h1 = alpha*(h2(x-1,y)+h2(x+1,y)-2h2(x,y)+ h2(x,y-1)+h2(x,y+1)-2h2(x,y)) - G \qquad (3)$$

We can let

$$f4 = (h2(x,y) - h1(x,y) + G)/\text{alpha} + 4*h2(x,y) =$$
$$h2(x-1,y) + h2(x+1,y) + h2(x,y-1) + h2(x,y+1). \quad (4)$$

f4 can also be viewed as the average of 4 times the center point. Using the gradually varied function, the average is $((h2-h1+G)/\text{alpha} + 4*h2(x,y))/4$. Based on the properties of gradually varied functions, the maximum difference is three. For example, h2_old(x-1),y) is bigger, then h2_new(x-1),y) should be smaller.

The iterating procedure is to update the center point after the first fit.

## 6. Real Data Processing and Application

This results use an algorithm to fit the initial data set using an individualized fit, see Fig 3. This algorithm is also made by the rough graduate varied surface fitting by scanning through the fitting array. There are many clear boundary lines in the images. In order to reduce error, our new algorithm will use more accurate formulas to calculate the derivatives (1):

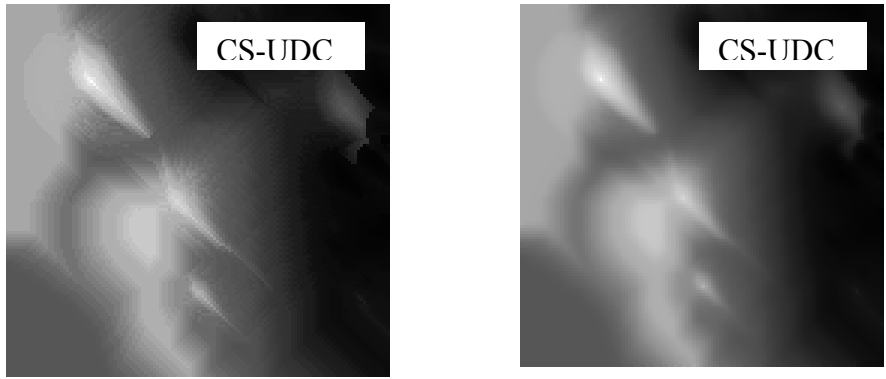

Fig. 3 Gradually varied functions for initial surface

Fig.4 shows the sequential surface fitting and the water flow equation updating. Individual surface fitting results are shown in (a), (b), and (c). Starting with fitted surface at each time the process will be faster to converge. It will not affect to the final result if there are enough iterations.

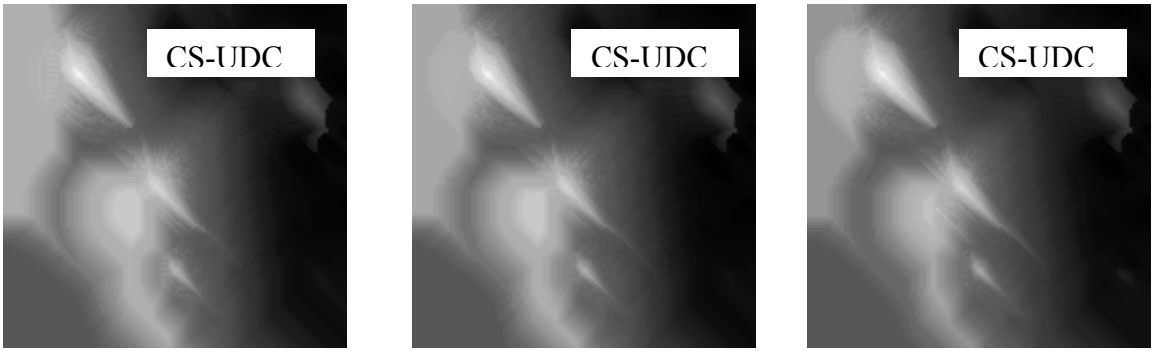

| (a) day1 | (b) day 30 | (c) day 50 |

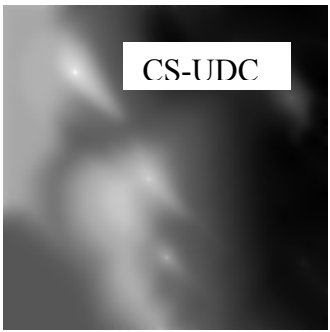 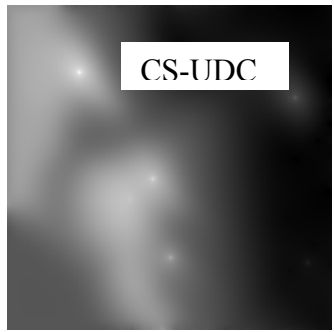 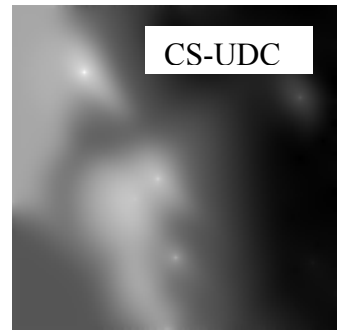

(d) Day 3 with flow equation iteration  (f) Day 30 with iteration  (g) Day 50 with iteration

Fig 4. VA Groundwater distribution calculated by gradually varied surfaces and the flow equation.

To test the correctness of the result, we have located the data points in the rectangle area in GoogleMap and Find latitude and Longitude. http://findlatitudeandlongitude.com/
We also tested the following fitting points that show the corrected correspondences.

Selected Points used in reconstruction
One can find the location at

| | | |
|---|---|---|
| 4.65 | 36.62074879 | -76.10938540 |
| 75.37 | 36.92515020 | -77.17746768 |
| 6.00 | 36.69104276 | -76.00948530 |
| 175.80 | 36.78431615 | -76.64328700 |
| 168.33 | 36.80403855 | -76.73495750 |
| 157.71 | 36.85931567 | -76.58634110 |
| 208.26 | 36.68320624 | -76.91329390 |
| 7.26 | 36.78737704 | -76.05153760 |

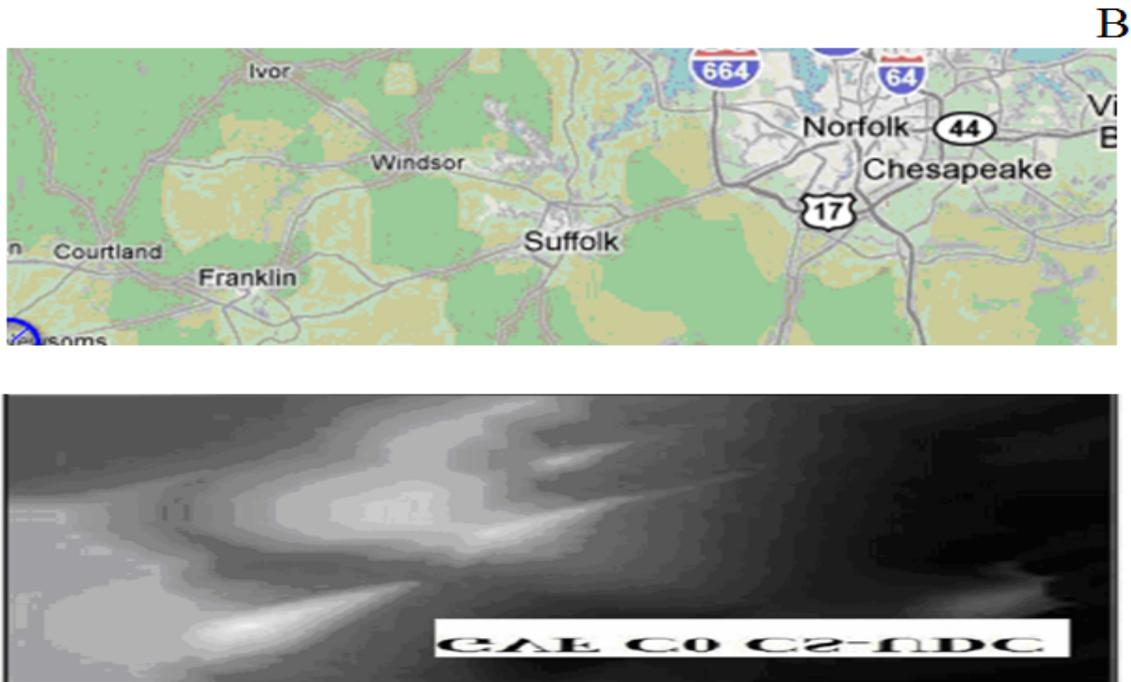

Fig. 5 The map and ground water data

Fig. 5 shows a good match found between the ground water data and the region geographical map. The brighter pixels mean the deeper distance from the surface. In mountain area, the groundwater level is lower in general. Some mismatches may be caused by not having enough sample data points (wells).

**Future work:** to ensure the accuracy of the calculation, we will add the finite element method to our research. We also want to use MODEFLOW to calculate local and small-region flow, and to use gradual variation to compute regional or global data.

*Acknowledgement:* This research was supported by USGS seed grants. The author expresses thanks to Dr. William Hare and members of UDC DC Water Resources Research Institute for their help. The author also thank Mr. Travis L. Branham for data collection. This paper is based on the report entitled: Li Chen, Gradual Variation Analysis for Groundwater Flow of DC, DC Water Resources Research Institute Final Report 2009 [24].

.

20. Travis L. Branham, Development of a Web-based Application to Geographically Plot Water Quality Data, UDC, 2008
21. G. Agnarsson and L. Chen, On the extension of vertex maps to graph homomorphisms, Discrete Mathematics, Vol 306, No 17, pp 2021-2030, Sept. 2006.
22. G.W. Pruist, B. H. Gilding and M. J. Peters, A comparison of different numerical methods for solving the forward problem in EEG and MEG, Physiol. Meas. 14, A1-A9, 1993.
23. Li Chen, Digital-Discrete Method For Smooth-Continuous Data Reconstruction, Submit to Capital Science 2010 of The Washington Academy of Sciences and its Affiliates, March 27 – 28, 2010.
24. Li Chen, Gradual Variation Analysis for Groundwater Flow of DC, DC Water Resources Research Institute Final Report 2009.

12